\documentclass[a4paper,12pt]{article}
\usepackage[utf8]{inputenc}

\usepackage{amsmath,graphicx}

\newtheorem{lemma}{Lemma}
\newtheorem{corollary}{Corollary}
\newtheorem{theorem}{Theorem}
\newtheorem{proposition}{Proposition}

\newenvironment{myp}{{\noindent \it \bf \sf Proof:\/}}{\hspace*{1em}\hfill  \frame{o}\vskip 0.1in}

\newcommand{\pf}{{\sf{pf}}}
\title{Wide enough Latin rectangles are perfect}
\author{N. Astromujoff and M. Matamalafootnote
\footnote{This work 
has been partially supported by program Basal-CMM (M.M.) and N\'ucleo Milenio Información y Coordinación en Redes ICM/FIC P30(N.A. and M.M.)}}
\author{Natacha Astromujoff$^1$ and Martín Matamala$^{1,2}$\\
\small (1) Departamento de Ingeniería Matemática\\
\small (2) Centro de Modelamiento Matemático (UMI 2807 CNRS),
\small \\Universidad de Chile
}

\begin{document}

\maketitle

\begin{abstract}
Given two integers $m$ and $n$ with $m\leq n$, 
a Latin rectangle of size $m\times n$  
is a bi-dimensional array with $m$ rows and $n$ columns
filled with symbols from an alphabet with $n$ symbols,
such that each row contains a permutation of the alphabet
and each column contains no repeated symbols.

Two rows $a$ and $b$ of a Latin rectangle $R$ 
define a permutation $R_{a,b}$ assigning the 
symbol $y$ to the symbol $x$ if they are in the same column, $x$ is in row $a$ and $y$ is in row $b$. 
A Latin rectangle $R$ is perfect is the permutation $R_{a,b}$
is cyclic, for each pair of rows $a$ and $b$.

We prove that for each integer $m$ and each large enough odd integer 
$n$ there is a perfect Latin rectangle $R$ of size $m\times n$.
It is a partial (asymptotic) answer to a well-known conjecture which says that 
the same property holds for each odd integer $m\leq n$.
\end{abstract}

\section{Preliminaries}

Let $m$ and $n$ be two integers with $m\leq n$.
A \emph{Latin rectangle} $R$ of size $m\times n$ is a matrix filled with symbols
from an alphabet $\Sigma_R$ of size $n$ such that each row of $R$ has a permutation of $\Sigma_R$ and each column contains $m$ distinct symbols.

We denote by $R(a,c)$ the symbol of $R$ in row $a$ and column $c$. 
Given two rows $a$ and $b$ of a Latin rectangle $R$ we define the permutation 
$R_{a,b}$ on $\Sigma_R$, by $R_{a,b}(x)=y$ if and only
if there is a column $c$ of $R$ such that $R(a,c)=x$ and $R(b,c)=y$. 
A pair of rows $(a,b)$ is \emph{perfect} if $R_{a,b}$ is a cyclic permutation.
We denote by $\pf(R)$ the number of perfect pairs of a Latin rectangle $R$. Then, $\pf(R)\leq \binom{m}{2}$, for a Latin rectangle with $m$ rows.
When, the latter inequality holds as equality we say that $R$ is a \emph{perfect} Latin rectangle. Let $A$ be the set of all integers $m$ such that there is a perfect Latin rectangle $S$ of size $m\times m$ (a Latin \emph{square} of \emph{order} $m$).

The study of perfect Latin squares is deeply related to \emph{perfect} one-factorizations of complete and complete bipartite graphs. In this latter context a one-factorization is \emph{perfect} if the graph induced by any pair of 1-factors
(which, in general, is a 2-factor) is a Hamiltonian cycle.

It is known that if a perfect one-factorization of the complete
graph $K_{n+1}$ exists, then also the complete bipartite graph $K_{n,n}$
has a perfect one-factorization (see \cite{wi2005} for a nice presentation).
It is also known that $K_{n,n}$ has 
a perfect one-factorization only if $n$ is odd or $n=2$. These two
results seem to be first proved in \cite{L1980}. More recently, it was established 
that $K_{n,n}$ has a perfect one factorization if and only $n\in A$
\cite{W1999}.

For all even $n\leq 52$ it is known that 
$K_n$ has a perfect one-factorization. 
Laufer \cite{L1980} suggests that Kotzig was who noticed that for each odd prime $p$, $K_{p+1}$ and $K_{2p}$ have perfect one-factorizations.
From these facts one gets that for each odd prime $p$, the graphs $K_{p,p}$ and $K_{2p-1,2p-1}$ have perfect one-factorizations. Later in \cite{BMW2002}, it was proved that also $K_{p^2,p^2}$ has a perfect one-factorization,
for each odd prime $p$. These results
imply that for  each odd prime $p$, perfect Latin rectangles of size $m\times n$ exist for each $n\in \{p,2p-1,p^2\}$ and each $m\leq n$.

It is conjectured that $K_n$ has a perfect one-factorization for each even $n$.
A weaker conjecture is that $K_{n,n}$ has a perfect one-factorization
for each odd $n$. This conjecture can be stated, equivalently, by saying that 
for each odd integer $n$ and each $m\leq n$, there is a perfect Latin rectangle
of size $m\times n$. The equivalence between both formulations is obvious as from each perfect Latin square of order $n$ we can obtain a perfect Latin rectangle of size $m\times n$, for each $m\leq n$.
\medskip 

\section {Our contribution}

In this work we prove that for each  integer $m$ and each odd $n$ large
enough, a perfect Latin rectangle exists of size $m\times n$.
More precisely, we prove the following.

\begin{theorem}\label{t:tech}
For each integer $m$ there is $n_0$ such that for all
$n\geq n_0$ there is a perfect Latin rectangle of size $m\times n$.
\end{theorem}

Let $B$ be the set of all integers $m$ satisfying Theorem \ref{t:tech}.
As from a perfect Latin rectangle of size $m\times n$ we can 
obtain perfect Latin rectangles of size $m'\times n$, for each
$m'\leq m$, in order to proof Theorem \ref{t:tech}, it is
enough to prove that $B$ is infinite. This shall prove this latter property 
by proving that the set $A$ is a subset of the set $B$.
In fact, the set $A$ is infinite as it contains the set of all primes.
Hence, if $A\subseteq B$, then $B$ is infinite as well.

The main technical contribution of this work is the following.

\begin{proposition}\label{l:widthextension}
Let $m\in A$. If a perfect Latin rectangle of size $m\times n$ exists, 
then also there is  a perfect Latin rectangle of size $m\times (n+m-1)$.
 \end{proposition}

The repeated application of Proposition \ref{l:widthextension}
implies that if $m\in A$ and a perfect Latin rectangle of size $m\times n$ exists, then also perfect Latin rectangles exist for each
$n'\geq n$ such that $n'\equiv n \mod (m-1)$. 
Therefore, thanks to Proposition \ref{l:widthextension} the proof of Theorem \ref{t:tech} is a consequence of the following property.
\begin{proposition}\label{modular}
For each $m\in A$ and each odd integer $i$ in $\{1,\ldots,m-2\}$, there is 
a perfect Latin rectangle of size $m\times n_i$, where $n_i$ is an odd integer 
such that $n_i \equiv i \mod (m-1)$ and $m\leq n_i$.
\end{proposition}

For $i=m-2$,  from Dirichlet's Theorem we known that there is a prime $p$ such that $p\equiv m-2 \mod (m-1)$ and $m\leq p$. As any prime belongs to $A$, it follows that there is a perfect Latin rectangle of size $m\times p$ and we can set $n_{m-2}=p$. Unfortunately, we cannot use Dirichlet's Theorem in all the remaining cases because it only applies when $i$ and $m-1$ are coprimes. 
However, we can use Proposition \ref{l:widthextension} and the case
$i=m-2$ in order to solve the remaining cases. In fact, we have the following result.

\begin{lemma}\label{l:general}
 Let $m$ be an odd integer  and let $r\in A$ such that $r\equiv m-2 \mod m-1$
 and $m\leq r$.
 Then, for each odd $i$ in $\{1,\ldots,m-2\}$
 there is a perfect Latin rectangle of size $m\times n_i$
 such that $n_i \equiv i \mod (m-1)$ and $m\leq n_i$.
 \end{lemma}

\begin{myp}
As $r\in A$ we can apply Proposition \ref{l:widthextension} to  $r$.
Then,  for each $j\geq 0$ we have that there is a perfect Latin rectangle of size $r\times (r+j(r-1))$. Since $m\leq r$ we get that for each $j\geq 0$, there is 
a perfect Latin rectangle of size $m\times (r+j(r-1))$. In order to finish 
the proof it is enough to find, for each odd $i\in \{1,\ldots,m-2\}$ a $j$ such that $m\leq r+j(r-1) \equiv i \mod (m-1)$. For each such $i$ we set $j=(m-2-i)/2$. 

By hypothesis we have that $r\equiv m-2 \mod m-1$.
Hence, for each odd $i$ in $\{1,\ldots, m-2\}$ we have that 
$$ r+(m-2-i)(r-1)/2\equiv m-2+(-2) (m-2-i)/2\equiv i \mod m-1.$$
 Therefore,   we obtain the conclusion by defining
$n_i:=r+(m-2-i)(r-1)/2$.

\end{myp}

We now give the proof of Proposition \ref{l:widthextension}.
 \begin{figure}\label{f}
\scriptsize
\begin{tabular}{c c c}
 \begin{tabular}{|c| c| c| c| c|}
 \hline
 \multicolumn{5}{|c|}{$R$}\\
 \hline 
 0 & 1 & 2 & {\bf 3} & 4 \\ \hline
4 &0 & 1 & {\bf 2} & 3  \\ \hline
3 & 4 & 0  & {\bf 1} & 2 \\ \hline
2 & 3 & 4 & {\bf 0} & 1 \\ \hline
 1 & 2 & 3 & {\bf 4} & 0 \\
  \hline
\end{tabular}
&
 \begin{tabular}{|c| c| c| c| c|}
 \hline
 \multicolumn{5}{|c|}{$S$}\\
 \hline 
 {\bf 5} & 6 & 7 & 8 & 9 \\ \hline
9 &{\bf 5} & 6 & 7 & 8  \\ \hline
8 & 9 & {\bf 5} & 6 & 7 \\ \hline
7 & 8 & 9 & {\bf 5} & 6 \\ \hline
 6 & 7 & 8 & 9 & {\bf 5} \\
  \hline
\end{tabular}
& $\to$
 \begin{tabular}{|c| c| c| c|| c| c| c| c| c|}
  \hline
 \multicolumn{9}{|c|}{$T$}\\
 \hline
  0 & 1 & 2  & 4 & {\bf 3} & 6 & 7 & 8 & 9 \\ \hline
4 &0 & 1  & 3  & 9 &{\bf 2} & 6 & 7 & 8  \\ \hline
3 & 4 & 0   & 2 & 8 & 9 & {\bf 1}  & 6 & 7 \\ \hline
2 & 3 & 4 & 1&7 & 8 & 9 & {\bf 0} & 6  \\ \hline
 1 & 2 & 3 & 0  &6 & 7 & 8 & 9 & {\bf 4}  \\ \hline
 \end{tabular}
\end{tabular}
\medskip 

\begin{tabular}{c c}
 \begin{tabular}{|c| c| c| c| c| c| c| c| c|}
  \hline
 \multicolumn{9}{|c|}{$T$}\\
 \hline
  {\bf 0} & 1 & 2  & 4 & {3} & 6 & 7 & 8 & 9 \\ \hline
{\bf 4} &0 & 1  & 3  & 9 &{2} & 6 & 7 & 8  \\ \hline
{\bf 3} & 4 & 0   & 2 & 8 & 9 & { 1}  & 6 & 7 \\ \hline
{\bf 2} & 3 & 4 & 1&7 & 8 & 9 & { 0} & 6  \\ \hline
 {\bf 1} & 2 & 3 & 0  &6 & 7 & 8 & 9 & { 4}  \\ \hline
 \end{tabular}
&    
 \begin{tabular}{|c| c| c| c| c|}
 \hline
 \multicolumn{5}{|c|}{$S'$}\\
 \hline 
 {10} & 11 & 12& 13 & {\bf 14} \\ \hline
{\bf 14} & 10 & 11 & 12 & 13  \\ \hline
13 & {\bf 14} & 10 & 11 & 12 \\ \hline
12 & 13 & {\bf 14} & 10  & 11 \\ \hline
 11 & 12 & 13 & {\bf 14} & 10 \\
  \hline
\end{tabular} $\to $
\end{tabular}
\medskip 

\ \  \begin{tabular}{|c| c| c| c| c| c| c| c|| c| c| c| c| c|}
  \hline
 \multicolumn{13}{|c|}{$T'$}\\
 \hline
 1 & 2  & 4 & {3} & 6 & 7 & 8 & 9 & 10 & 11& 12& 13 & {\bf 0} \\ \hline
0 & 1  & 3  & 9 &{2} & 6 & 7 & 8  &{\bf 4} & 10 & 11& 12& 13 \\ \hline
 4 & 0   & 2 & 8 & 9 & { 1}  & 6 & 7  & 13& {\bf3}& 10& 11&12 \\ \hline
 3 & 4 & 1&7 & 8 & 9 & { 0} & 6   & 12 & 13& {\bf 2} &10 &11 \\ \hline
 2 & 3 & 0  &6 & 7 & 8 & 9 & { 4}   & 11&12 &13 & {\bf 1} & 10 \\ \hline
 \end{tabular}

\caption{Two applications of Proposition \ref{l:widthextension}.
In the first application we use two copies of a perfect Latin square of 
order $m=n=5$ as $R$ and $S$. We delete the fourth column of $R$ and the symbol $\bf 5$ of $S$.
In the second application we use the resulting Latin rectangle $T$ with 
another copy of the Latin square $S$. We delete the first column of $T$ 
and the symbol $\bf 14$ of $S'$.}
\end{figure}

\begin{myp} Let $R$ be a perfect Latin rectangle of size $m\times n$
and let $S$ be a perfect Latin square of order $m$. We build a perfect
Latin rectangle $T$ of size $m\times (n+m-1)$. Let us assume that the rows of $R$ and $S$ are indexed by the same set of indices $I$, and that 
their columns are indexed by disjoint sets $J_R$ and $J_S$, respectively.
We also assume that $\Sigma_R$ and $\Sigma_S$ are disjoint.

Let $c\in J_R$ arbitrary and let $R'$ be the matrix of size $m\times (n-1)$ obtained from $R$ by deleting column $c$. Then, the symbol $R(a,c)$ is missing
in $R'$ at row $a$, for each $a$ in $ I$.
Let $s\in \Sigma_S$ arbitrary and let $S'$ be the matrix obtained from $S$ by replacing symbol $s$ at row $a$ by $R(a,c)$. Finally, let $T$ be the matrix of size $m\times (n+m-1)$ obtained
by concatenating $R'$ with $S'$ (see Figure \ref{f} for an example).
Then, the set of symbols of $T$ is $\Sigma_T=\Sigma_R \cup (\Sigma_S\setminus \{s\})$.
The set of columns of $T$ is indexed by the set $J_R\setminus \{c\} \cup J_S$.

Let $S(a)$ be the column of $S$ such that $S(a,S(a))=s$.
Formally, $T$ is given by 
$T(a,d)=R(a,d)$ if $d\in J_R\setminus \{c\}$,
$T(a,d)=S(a,d)$ if $d\in J_S\setminus\{S(a)\}$,  and 
$T(a,S(a))=R(a,c)$, for each $a$ in $I$.

Since $\Sigma_R$ and $\Sigma_S$ are pairwise disjoint it is clear that 
$T$ is a Latin rectangle.

Let $a$ and $b$ be two rows in $I$. 
We prove that the permutation $T_{a,b}$ is cyclic by proving that 
$T_{a,b}^t(R(b,c))\neq R(b,c)$ for each $t<n+m-1$ and that 
$T_{a,b}^{n+m-1}(R(b,c))=R(b,c)$.

From the definition of $T$ and since $R$ is a perfect Latin rectangle we have that for each $t<n-1$, $T_{a,b}^t(R(b,c))=R_{a,b}^t(R(b,c))$ and that $T_{a,b}^{n-1}(R(b,c))=R(a,c)$. In fact, the iterated application of $T_{a,b}$ to $R(b,c)$ has the same effect  that has the application of $R_{a,b}$, untill the symbol $R(a,c)$ is reached.

The symbol $R(a,c)$ is in row $a$ and column $S(a)$
of $T$ hence $T_{a,b}^{n}(R(b,c))=T_{a,b}(R(a,c))=S(b,S(a))$.
Since $S$ is a perfect Latin square, the behavior of $T_{a,b}$ in the 
next $m-2$ iterations is the same as that of $S_{a,b}$. Therefore, after $n+m-2$ iterations
of $T_{a,b}$ over $R(b,c)$ we obtain the symbol $S_{a,b}^{m-2}(S(b,S(a)))=S(a,S(b))$. From the definition of $T$ we have that $T(S(a,S(b)))=R(b,c)$ and then after $n+m-1$ iterations of $T_{a,b}$ 
we come back to symbol $R(b,c)$.

\end{myp}

We now discuss how to get upper bounds for the value of $n_0$ appearing
in Theorem \ref{t:tech}. Given an odd integer $m$ and an odd integer $i$ in $\{1,\ldots,m-1\}$, let $\theta(m,i)$ be the smallest integer $k$ such that 
$k\equiv i \mod (m-1)$, $m\leq k$  and there is a perfect Latin rectangle of
size $m\times k$. From this definition for each $m\in A$ we have that $\theta(m,1)=m$.
For any integer $m\geq 2$, let  $\theta(m)$ denote the smallest integer $k$ such that $m\leq k$ and for each $n\geq k$ there is a perfect Latin rectangle of size $m\times n$. Then $\theta(m)$ is the smallest $n_0$ for which the statement of 
Theorem \ref{t:tech} is valid.
From Proposition \ref{l:widthextension},  it is easy to see that  for each $m\in A$ we have that  $\theta(m)\leq \max\{\theta(m,i): i \mbox{ odd in }\{1,\ldots,m-2\}\}$.

For small values of odd $m$, we know that $\theta(m)=m$.
In fact, for each odd $m\leq 27$, every integer in $\{m,m+2,\ldots,2m-3\}$ belongs to $A$. Hence, for $m\leq 27$ and each odd $i$ in $\{1,\ldots,m-2\}$ we have that $\theta(m,i)=m-1+i$. Thus, $\theta(m)=m$ for each odd $m\leq 27$.

More generally, in \cite{L1944a,L1944b} it was proved a quantitative version of Dirichlet's Theorem. In terms of our discussion, it says that 
there is an universal constant $L$ such that $\theta(m,m-2)\leq (m-1)^{L}$, for each odd $m$. The best known value for $L$ is $5.2$ \cite{X2011}.  
It is also known that this bound can be improved to $m^2(\log m)^2$ if
the so called Strong Riemann Hypothesis holds. From the proof of Lemma \ref{l:general}  both bounds for $\theta(m,m-2)$ transfer directly to bounds for $\theta (m)$: $m^{6.2}$ and $m^3(\log m)^2$, respectively,  when $m\in A$ since in this case we can apply
Proposition  \ref{l:widthextension}. When $m$ is an arbitrary integer we can
use a classical result of Chebyshev \cite{C1852} saying that there is a prime $m'$ with $m\leq m'\leq 2m$. As each prime $m'$ belongs to $A$, from previous analysis we get that $\theta(m)\leq \theta(m')\leq (m')^{6.2}\leq 74m^{6.2}$, for each $m$. We summarize this discussion in the following corollary.

\begin{corollary}\label{c:quan}
 For every integer $m$ and every odd $n\geq 74 \lfloor m^{6.2}\rfloor$
 there is a perfect Latin rectangle of size $m\times n$.
\end{corollary}

\end{document}